\documentclass[aps, preprint,unsortedaddress,
numerical
]{revtex4-1}

\usepackage{hyperref}
\usepackage{bm}

\newtheorem{theorem}{Theorem}[section]
\newtheorem{lemma}[theorem]{Lemma}

\usepackage{amsmath,amssymb}

\newcommand{\qed}{\hfill$\square$}

\begin{document}

%\preprint{}

\title[]{Affine Killing vector fields on homogeneous surfaces with torsion}

\author{D.~D'Ascanio}
\homepage{dascanio@fisica.unlp.edu.ar}
\affiliation{%
Instituto de F\'isica La Plata, CONICET and\\ Universidad Nacional de La Plata, CC 67 (1900) La Plata, Argentina.}

\author{P.~B.~Gilkey}
\homepage{gilkey@uoregon.edu}
\affiliation{%
Mathematics Department, University of Oregon, Eugene OR 97403 USA.}

\author{P.~Pisani}
 \homepage{pisani@fisica.unlp.edu.ar}
\affiliation{%
Instituto de F\'isica La Plata, CONICET and\\ Universidad Nacional de La Plata, CC 67 (1900) La Plata, Argentina.}

\date{\today}

\begin{abstract}
Many extensions of General Relativity are based on considering metric and affine structures as independent properties of spacetime. 
	This leads to the possibility of introducing torsion as an independent degree of freedom. In this article we examine the effects of torsion 
	on the affine Killing vectors of two-dimensional manifolds. We give a complete description of the Lie algebras of affine Killing vector 
	fields on homogeneous surfaces. This can be used in the search of non-metrizable surfaces of interest.
\end{abstract}

\maketitle

\section{Introduction}

General Relativity is at present the most successful description of the gravitational interaction. However, many open questions settled by the present status of astrophysical observations 
motivate the search for modified formulations of this theory. In addition, black hole physics 
and early universe models require a framework compatible with quantum mechanics. For these 
reasons, General Relativity in its present form is not considered as an ultimate description of gravity 
and different generalizations are currently under study.

One approach to this reformulation is based on a reexamination of the canonical degrees of freedom of the theory. 
Constructing an invariant action requires a metric and an affine connection, both describing different
geometric properties of spacetime. In standard General Relativity, only the metric is a fundamental field whereas 
the affine structure is given by the Levi-Civita connection. However,  from the mathematical point of view,
the Riemannian and affine structures need not be
related; the connection is an independent degree of freedom locally given by non-metrizable Christoffel 
symbols \cite{Zanelli:2005sa}. In a general setting, the difference between the Christoffel symbols and 
those derived from the Levi-Civita connection is given by the non-metricity tensor and the torsion tensor \cite{Hehl:1994ue}.

Geometries with non-vanishing non-metricity have attracted renewed attention aimed at exploring the coupling to 
matter \cite{Iosifidis:2018jwu} as well as the geometric properties of new spacetime configurations (see e.g. \cite{Klemm:2018bil} and references therein). On the other 
hand, gravity theories with non-vanishing torsion have been the subject of extensive study: There are models of the early 
universe in which torsion has a fundamental role as an alternative to inflation \cite{Trautman:1973wy,Poplawski16}. 
The propagation of quantum fields on a spacetime with torsion has been analyzed in \cite{Shapiro:2001rz}, and an 
example of the torsion field as a propagating degree of freedom can be found in \cite{Blagojevic14}. Other models with non-metricity and torsion waves have been studied in \cite{Babourova:2018crn} and \cite{Babourova:1998ct}. A one-loop effective action in terms of the connection has been analyzed in \cite{YuBaurov:2018pyj}. Apart from the 
theoretical interest in non-Riemannian models of gravity, there are a variety of experiments which have been 
designed to measure the effects of torsion; for a quite comprehensive account we refer to \cite{Hammond:2002rm}.
There is also a substantial body of literature in the purely mathematical setting (we refer to \cite{BG18,BK17,GM17,KM16,M18,Y17}, to cite just a few representative examples).

In this context we consider it useful to pursue the analysis of purely affine properties 
without regard to any possible Riemannian structure. The purpose of this work is to examine 
how the torsion impacts the geometry of a surface; we shall focus our attention on describing the effect of torsion on the  associated Lie algebra of affine Killing vectors. As we do not use field equations, our results are model independent.

To ensure that the Lie algebra of affine Killing vectors is sufficiently rich, we shall assume that the surface in question
is locally homogeneous; a complete classification of such local geometries is available \cite{AMK,Opozda}.
In addition, we believe that our study gives insight into the analysis of 
three- and higher-dimensional manifolds. This notwithstanding, theories of gravity with torsion in two dimensions constitute an active area on their own---for a review of its motivations and development see \cite{Grumiller:2002nm,Obukhov:1997uc,Katanaev:1986qm}.

In the present paper we shall assume the underlying manifold in question is simply connected
to facilitate the passage from local to global questions. The Lie algebra $\mathfrak{K}$
of affine Killing vector fields has played an important role in 
the study of surfaces which are torsion free; in this paper, we examine the relationship between the torsion and $\mathfrak{K}$.
We say that a Lie sub-algebra $\mathfrak{K}_0$ of $\mathfrak{K}$ is {\it effective} if given any point $P$ of the
underlying manifold, there exist $X_i\in\mathfrak{K}_0$ so that $\{X_1(P),X_2(P)\}$ are linearly independent. Since the underlying structure is assumed locally homogeneous and simply connected, $\mathfrak{K}$ is effective 
(see Hall~\cite{Hall} or Nomizu~\cite{Nomizu}). We refer to \cite{BiG17,V16} for recent examples where affine Killing vector fields
have played an important role in the analysis.

We first present the fundamental definitions and properties of affine manifolds introducing torsion and the space of affine
Killing vector fields. We focus on homogeneous affine surfaces and recall known results concerning their classification, 
in particular those related to affine Killing vector fields. We state the main result of the paper, namely the description of locally homogeneous affine surfaces in terms of the algebra of their affine Killing vectors.

\subsection{Affine surfaces, Christoffel symbols, and the torsion tensor}
An affine surface  is a pair $\mathcal{M}=(M,\nabla)$ where $M$ is a smooth surface and $\nabla$ is a connection on the tangent bundle of $M$. In contrast to the notation adopted by some authors, we emphasize that we permit $\nabla$ to have torsion. Let $\partial_k:=\frac\partial{\partial x^k}$
in some system of local coordinates $(x^1,x^2)$ on $M$.
We sum over repeated indices to express
$\nabla_i\partial_j = \Gamma_{ij}^k \partial_k$; the connection is determined by  
the Christoffel symbols $\Gamma_{ij}^k$. For two vectors $X,Y$, let $T(X,Y):=\nabla_XY-\nabla_YX-[X,Y]$ be
the {\itshape torsion tensor}; the components of the torsion tensor are expressed by
$$
T=(dx^i\wedge dx^j)\otimes(\Gamma_{ij}^k-\Gamma_{ji}^k)\partial_k
=(dx^1\wedge dx^2)\otimes {4T^i\partial_i}\text{ for }T^i:=\tfrac12(\Gamma_{12}^i-\Gamma_{21}^i)\,.
$$
We say $\mathcal{M}$ is {\itshape torsion free} if $T=0$, i.e. if
$\Gamma_{12}^k=\Gamma_{21}^k$ for $k=1,2$. There is a canonically associated torsion free connection
${}^0\nabla:=\nabla-T$ with Christoffel symbols 
$${}^0\Gamma_{ij}^k:=\textstyle\frac12(\Gamma_{ij}^k+\Gamma_{ji}^k)\,.
$$
The connection ${}^0\nabla$ is in a certain sense the symmetric part of the connection $\nabla$ and the torsion $T$
is the anti-symmetric part. We let ${}^0\!\mathcal{M}:=(M,{}^0\nabla)$. 
Conversely, given an affine surface without torsion ${}^0\!\mathcal{M}$ and a torsion tensor $T=(dx^1\wedge dx^2)\otimes({4T^i\partial_i})$, we can perturb ${}^0\!\mathcal{M}$ to
define a surface ${}^T\!\mathcal{M}$ with torsion $T$ by setting ${}^T\nabla={}^0\nabla+T$; the resulting Christoffel
symbols are given by setting:
\begin{equation}\label{E1.a}
\begin{array}{llll}
{}^T\Gamma_{11}^1={}^0\Gamma_{11}^1,&{}^T\Gamma_{11}^2={}^0\Gamma_{11}^2,&{}^T\Gamma_{22}^1={}^0\Gamma_{22}^1,\\[0.03in]
{}^T\Gamma_{22}^2={}^0\Gamma_{22}^2,&
{}^T\Gamma_{12}^1={}^0\Gamma_{12}^1+T^1,&{}^T\Gamma_{12}^2={}^0\Gamma_{12}^2+T^2,\\
{}^T\Gamma_{21}^1={}^0\Gamma_{12}^1-T^1,&{}^T\Gamma_{21}^2={}^0\Gamma_{12}^2-T^2.
\end{array}\end{equation}
These constructions are independent of the particular coordinate system chosen. 

\subsection{Affine Killing vector fields} 
Let $\mathcal{M}$ be an affine surface. A smooth vector field $X=v^1\partial_1+v^2\partial_2=(v^1,v^2)$ on an affine surface is said to be an 
{\itshape affine Killing vector field} if the Lie derivative 
of the connection with respect to the vector field $X$ vanishes or, 
equivalently (see Kobayashi and Nomizu~\cite[Chapter VI]{KN63}), if the 8 affine Killing equations for $1\le i,j,k\le2$ are satisfied
\begin{align}\label{E1.b}
K_{ij}^k:\quad 0={} \frac{\partial^2 v^k}{\partial x^i \partial x^j} + v^l \frac{\partial \Gamma_{ij}^k}{\partial x^l} 
- \Gamma_{ij}^l\frac{\partial v^k}{\partial x^l} + \Gamma_{il}^k\frac{\partial v^l}{\partial x^j}
+ \Gamma_{lj}^k\frac{\partial v^l}{\partial x^i} \,.
\end{align}
The affine Killing equations form an over determined elliptic system of second-order partial differential equations.
The Lie bracket makes the linear space $\mathfrak{K}(\mathcal{M})$ of affine Killing vector fields into a Lie algebra of dimension
at most $6$ since an affine Killing vector field is determined by $X(0)$ and $\nabla X(0)$.

\subsection{Homogeneous affine surfaces} We say that a diffeomorphism from one affine surface to another is an
{\it affine map} if it intertwines the two associated connections.  
We say that an affine manifold $\mathcal{M}$ is {\itshape affine homogeneous}
if the Lie group of affine diffeomorphisms of
$\mathcal{M}$ acts transitively; the corresponding local notion is defined similarly. 
To pass between local and global results, we shall assume henceforth that the underlying manifold $M$ is simply connected and
locally affine homogeneous. In this setting, every affine Killing vector field which is locally defined extends to a
globally defined affine Killing vector field.

Opozda~\cite{Opozda} classified the locally homogeneous affine surfaces without torsion;
this classification was later extended to the case of surfaces with torsion by Arias-Marco and Kowalski~\cite{AMK}. We summarize their
result as follows.

\begin{theorem}\label{thm:opozda}
	If $\mathcal{M}$ is a locally homogeneous affine surface, possibly with torsion, then at least one of the following possibilities holds.
	\begin{enumerate}
		\item There exists a coordinate atlas for $M$ so that the Christoffel symbols of $\nabla$ are constant;
		$\mathcal{M}$ is said to be Type~$\mathcal{A}$.
		\item  There exists a coordinate atlas for $M$ so that the Christoffel symbols have the form 
		$\Gamma_{ij}^k = (x^1)^{-1} A_{ij}^k$, with $A_{ij}^k$ constant; $\mathcal{M}$ is said to be Type~$\mathcal{B}$.
		\item There exists a coordinate atlas for $M$ such that $\nabla$ is isomorphic to
		the Levi-Civita connection of the round sphere.
	\end{enumerate} 
\end{theorem}

We say that $\mathcal{M}=(\mathbb{R}^2,\nabla)$ is a {\it Type~$\mathcal{A}$ model} if the Christoffel symbols
$\Gamma_{ij}^k$ are constant.
If we identify $\mathbb{R}^2$ with the group of translations, then $\nabla$ is a Type~$\mathcal{A}$ model if
and only if $\nabla$ is left invariant.
We can describe Type $\mathcal{A}$ models in terms of the algebra of translations in the plane. Let $\mathfrak{K}_{\mathcal{A}}:=\operatorname{span}\{\partial_1,\partial_2\}$. Then $\mathcal{M}=(\mathbb{R}^2,\nabla)$
is a Type~$\mathcal{A}$ model if and only if $\mathfrak{K}_{\mathcal{A}}\subset\mathfrak{K}(\mathcal{M})$.
We say that $\mathcal{N}=(\mathbb{R}^+\times\mathbb{R})$ is a
{\it Type~$\mathcal{B}$ model} if $\Gamma_{ij}^k=(x^1)^{-1}A_{ij}^k$ for $A_{ij}^k$ constant. We identify 
$\mathbb{R}^+\times\mathbb{R}$
with the $ax+b$ group under the action $(x^1,x^2)\rightarrow(ax^1,ax^2+b)$; 
$(\mathbb{R}^+\times\mathbb{R},\nabla)$ is a Type~$\mathcal{B}$ model if and only if
$\nabla$ is left invariant under the natural action of the $ax+b$ group. 
By Theorem~\ref{thm:opozda}, any locally homogeneous
surface geometry is locally isomorphic to either a Type~$\mathcal{A}$ model, a Type~$\mathcal{B}$ model, or the round 2-sphere.
We remark that there are geometries which admit both Type~$\mathcal{A}$ and Type~$\mathcal{B}$ structures. We also
note that there are simply connected geometries with a Type~$\mathcal{A}$ structure which are not affine
equivalent to any open subset of a Type~$\mathcal{A}$ model; more than one coordinate system is required for such
geometries.

\subsection{The algebra of affine Killing vector fields for homogeneous surfaces}

Let $\mathcal{M}$ be a simply connected locally homogeneous affine surface. Fix a basepoint of $\mathcal{M}$.
Define the following Lie algebra structures on $\mathbb{R}^2$ and $\mathbb{R}^3$ by the nonzero brackets:
$$\begin{array}{llll}
\mathfrak{K}_{\mathcal{A}}:&[e_1,e_2]=0,\\[0.03in]
\mathfrak{K}_{\mathcal{B}}:&[e_1,e_2]=e_1,\\[0.03in]
\mathfrak{so}(3):&[e_1,e_2]=e_3,&[e_2,e_3]=e_1,&[e_3,e_1]=e_2,\\[0.03in]
\mathfrak{so}(2,1):&[e_1,e_2]=e_1,&[e_2,e_3]=e_3,&[e_3,e_1]=2e_2.
\end{array}$$
As already noted, $\mathfrak{K}_\mathcal{A}$ is the algebra of translations in the plane and
$\mathfrak{K}_\mathcal{B}$ is the algebra of horizontal translations and dilatations in the upper half-plane.
Following the notation of Patera et al.~\cite{Patera}, we define the following 4-dimensional Lie algebras by specifying their non-zero brackets: 
$$\begin{array}{lllll}
A_{4,9}^0:&[e_2,e_3]=e_1,&[e_1,e_4]=e_1,&[e_2,e_4]=e_2,\\[0.03in]
A_{4,12}:&[e_1,e_3]=e_1,&[e_2,e_3]=e_2,&[e_1,e_4]=-e_2,&[e_2,e_4]=e_1.
\end{array}$$
Let $\mathbb{A}_6$ be the 6-dimensional Lie algebra of the full affine group.

Recently Brozos-V\'azquez et al.~\cite{BVGRG19} gave a quite different proof Theorem~\ref{thm:opozda} by examining the affine Killing
equations directly. Their result, from which Theorem~\ref{thm:opozda} follows, may be stated as follows.

\begin{lemma}
	Let $\mathcal{M}=(M,\nabla)$ be locally homogeneous and simply connected. 
	\begin{enumerate}
		\item There is an effective Lie subalgebra $\tilde{\mathfrak{K}}$ of $\mathfrak{K}(\mathcal{M})$
		which is isomorphic to $\mathfrak{K}_{\mathcal{A}}$, $\mathfrak{K}_{\mathcal{B}}$, or $\mathfrak{so}(3)$.
		\item If $\tilde{\mathfrak{K}}\approx\mathfrak{K}_{\mathcal{A}}$, then
		there is a coordinate atlas so that $\Gamma_{ij}{}^k$ are constant.
		\item If $\tilde{\mathfrak{K}}\approx\mathfrak{K}_{\mathcal{B}}$, then
		there is a coordinate atlas so that
		$\Gamma_{ij}{}^k=(x^1)^{-1}A_{ij}{}^k$ for constant $A_{ij}{}^k$.
		\item If $\tilde{\mathfrak{K}}\approx\mathfrak{so}(3)$, then
		there is a coordinate atlas where
		$\nabla$ is the Levi-Civita connection defined by the metric of the round sphere.
\end{enumerate}\end{lemma}

In this paper we will complete their analysis. Our main result is the following; it is implicit in the computations
of Arias-Marco and Kowalski~\cite{AMK} but is not stated in this fashion there; our approach is quite different from theirs.
\begin{theorem}\label{T^1.3}
	Let $\mathcal{M}$ be a locally homogeneous simply connected affine surface with torsion.
	\begin{enumerate}
		\item Suppose $\mathcal{M}$ contains an effective Lie subalgebra which is isomorphic to $\mathfrak{K}_{\mathcal{A}}$.
		Then $\mathfrak{K}(\mathcal{M})$ is isomorphic to $\mathfrak{K}_{\mathcal{A}}$, to
		$\mathfrak{K}_{\mathcal{B}}\oplus\mathfrak{K}_{\mathcal{B}}$, to $A_{4,9}^0$, or to $A_{4,12}$.
		\item Suppose $\mathcal{M}$ contains an effective Lie subalgebra which is isomorphic to $\mathfrak{K}_{\mathcal{B}}$.
		Then $\mathfrak{K}(\mathcal{M})$ is isomorphic to $\mathfrak{K}_{\mathcal{B}}$, to
		$\mathfrak{K}_{\mathcal{B}}\oplus\mathfrak{K}_{\mathcal{B}}$, to $A_{4,9}^0$, or to $\mathfrak{so}(2,1)$.
		\item Suppose $\mathcal{M}$ contains an effective Lie subalgebra which is isomorphic to $\mathfrak{so}(3)$. Then
		$\mathcal{M}$ is without torsion and modeled on the round sphere.
\end{enumerate}\end{theorem}

\subsection{Outline of the paper}

The remainder of this paper is devoted to the proof of Theorem~\ref{T^1.3}. We begin in Section~\ref{S2} by establishing
the following useful observation. 

\begin{lemma}\label{L1.4}
	Let $\mathcal{M}$ be an affine surface and let ${}^0\!\mathcal{M}$ be the associated surface without torsion.
	Then $\mathfrak{K}(\mathcal{M}){{}\subseteq{}}\mathfrak{K}({}^0\!\mathcal{M})$. 
\end{lemma}

Brozos-V\'azquez et al. \cite{BVGRG18} and Gilkey and Valle-Regueiro~\cite{GVR} have classified, up to
linear equivalence, all the Type~$\mathcal{A}$ and Type~$\mathcal{B}$ models without torsion where
$\dim\{\mathfrak{K}\}>2$. Given an arbitrary model ${}^T\!\mathcal{M}$ of Type~$\mathcal{A}$ or Type~$\mathcal{B}$ with torsion, 
we pass to the associated torsion free model ${}^0\!\mathcal{M}$ and write down a basis for
$\mathfrak{K}({}^0\!\mathcal{M})$. We then examine the affine Killing equations to determine which affine Killing vector fields on ${}^0\!\mathcal{M}$ are affine Killing vector fields for $\mathcal{M}$. 
This then provides a classification of all the Type~$\mathcal{A}$ and Type~$\mathcal{B}$ models
with torsion where $\dim\{\mathfrak{K}(\mathcal{M})\}>2$, which is of interest in
its own right. Once this classification has been performed, we analyze the resulting Lie algebras
to complete the proof of Theorem~\ref{T^1.3}. This analysis is performed in Section~\ref{S3} in the Type~$\mathcal{A}$
setting and in Section~\ref{S4} in the Type~$\mathcal{B}$ setting. 
The original analysis of Brozos-V\'azquez et al.~\cite{BVGRG18}
ignored the flat geometries as being uninteresting as they are in the torsion free setting. But once torsion is added, it is necessary
to include these geometries as the flat geometries give rise to non-trivial geometries with torsion and for this the additional analysis of
Gilkey and Valle-Regueiro~\cite{GVR} is required.

%%%%%%%%%%%%%%%%%%%%%%%%%%%%%%%%%%%%%%%%%%%%%%%%%%

\section{Affine Killing equations in the presence of torsion}\label{S2}

In this section we give a proof of Lemma~\ref {L1.4}. We also give two examples which help to understand the role of torsion in the affine Killing algebra.

Let $(v^1, v^2)\in\mathfrak{K}({}^T\!\mathcal{M})$. Call ${}^T\!K_{ij}^k$ the r.h.s of Equation~(\ref{E1.b})
when the Christoffel symbols involve a torsion $T$. The corresponding equation for the symmetrized part of the Christoffel symbols is ${}^0K_{ij}^k$. A direct computation gives
$$
{}^T\!K_{ij}^k = {}^0K_{ij}^k + v^l\frac{\partial T_{ij}^k}{\partial x^l} - T_{ij}^l \frac{\partial v^k}{\partial x^l} + T_{il}^k \frac{\partial v^l}{\partial x^j} + T_{lj}^k \frac{\partial v^l}{\partial x^i} = 0\,.
$$
Taking $i=j$ in the last expression gives $^T\!K_{ii}^k = {}^0K_{ii}^k=0$. 
To obtain a similar result for the non-diagonal elements consider the equations
\begin{eqnarray*}
	&&{}^T\!K_{ij}^k = {}^0K_{ij}^k + v^l\frac{\partial T_{ij}^k}{\partial x^l} - T_{ij}^l \frac{\partial v^k}{\partial x^l} + T_{il}^k \frac{\partial v^l}{\partial x^j} + T_{lj}^k \frac{\partial v^l}{\partial x^i} = 0\,,\\
	&&
	{}^T\!K_{ji}^k = {}^0K_{ji}^k + v^l\frac{\partial T_{ji}^k}{\partial x^l} - T_{ji}^l \frac{\partial v^k}{\partial x^l} + T_{jl}^k \frac{\partial v^l}{\partial x^i} + T_{li}^k \frac{\partial v^l}{\partial x^j} = 0\,.
\end{eqnarray*}
Adding these we have $0={}^T\!K_{ij}^k + {}^T\!K_{ji}^k = {}^0K_{ij}^k + {}^0K_{ji}^k$. 
Since the Christoffel symbols for $T=0$ are symmetric, ${}^0K_{ij}^k = {}^0K_{ji}^k =0$. 
Lemma~\ref{L1.4} follows.
\qed

\subsection*{Example II.A}
Let $\mathcal{M}^4_1$ be the Type~$\mathcal{A}$ surface without torsion defined by the Christoffel symbols 
$\Gamma_{11}^1=-1$, $\Gamma_{12}^1 = 1$, $\Gamma_{22}^1=0$, $\Gamma_{11}^2 = 0$, 
$\Gamma_{12}^2=0$ and $\Gamma_{22}^2 = 2$. Let $0\ne T=(T^1,T^2)\in\mathbb{R}^2$. We shall see presently that 
$\dim\{\mathfrak{K}(\mathcal{M})\}=4$, that $\dim\{\mathfrak{K}({}^T\!\mathcal{M})\}=4$ if $T^2=0$. This shows that the equality in Lemma \ref{L1.4} can hold.

\subsection*{Example II.B}
Given a torsion free manifold which is locally homogeneous, the perturbed manifold  need not be homogeneous.
Consider the type $\mathcal{A}$ surface $\mathcal{M}^6_1$ defined by the Christoffel symbols 
$\Gamma_{11}^1=1$, $\Gamma_{12}^1 = 0$, $\Gamma_{22}^1=0$, $\Gamma_{11}^2 = 0$, 
$\Gamma_{12}^2=1$ and $\Gamma_{22}^2 = 0$, with $\dim\{\mathfrak{K}(\mathcal{M}^6_1)\}=6$. 
Perturb it by adding a type $\mathcal{B}$ torsion $T$ where $T^{1}=0$ and $T^{2} = t^{2}/x^1\neq0$. The resulting 
structure has  
$$
{\mathfrak{K}({}^T\!\mathcal{M}^6_1)=\mathrm{span}\{ \partial_2, x^2\partial_2, e^{-x^1}\partial_2 \}\,.}
$$
This algebra has no effective subalgebras for all $(x^1,x^2)$ and hence the surface is not homogeneous. 
Now perturb it by adding a type $\mathcal{B}$ torsion $T$ with $T^{1}=t^{1}/x^1\neq0$ and $T^{2} = 0$. The resulting 
structure has  $\mathfrak{K}({}^T\!\mathcal{M}^6_1)=\mathrm{span}\{ \partial_2 \}$.
This example shows that the addition of torsion to a homogeneous, torsion free surface does not necessarily give a homogeneous surface.

%%%%%%%%%%%%%%%%%%%%%%%%%%%%%%%%%%%%%%%%%%%%%%%%%%

\section{Type~$\mathcal{A}$ surfaces with torsion}\label{S3}

In this section we obtain the spaces of affine Killing vector fields for Type $\mathcal{A}$ models. This gives the algebras of Theorem~\ref{T^1.3}~(1).

Parametrize the set of Type~$\mathcal{A}$ models by setting
$\mathcal{M}(\vec\xi):=(\mathbb{R}^2,\nabla^{\mathcal{A}}(\vec\xi))$ for $\xi\in\mathbb{R}^8$
where the Christoffel symbols of $\nabla^{\mathcal{A}}(\vec\xi)$ are given by:
$$\begin{array}{llll}
\Gamma_{11}{}^1=\xi_1,&\Gamma_{11}{}^2=\xi_2,&\Gamma_{12}{}^1=\xi_3,&\Gamma_{12}{}^2=\xi_4,\\
\Gamma_{21}{}^1=\xi_5,&\Gamma_{21}{}^2=\xi_6,&\Gamma_{22}{}^1=\xi_7,&\Gamma_{22}{}^2=\xi_8.
\end{array}$$
The torsion free models $\mathcal{M}(\vec\xi)$ form
a 6-dimensional subspace where $\xi_3=\xi_5$ and $\xi_4=\xi_6$. The general linear group $\operatorname{GL}(2,\mathbb{R})$
acts on the space of Type~$\mathcal{A}$ models by change of basis and defines thereby a linear representation
of  $\operatorname{GL}(2,\mathbb{R})$ on $\mathbb{R}^8$. We say that two Type~$\mathcal{A}$ models are {\it linearly
	equivalent} if they lie in the same orbit of this representation.
The works \cite{BVGRG18,GVR} mentioned
previously classifies all Type~$\mathcal{A}$ torsion free
models up to linear equivalence. We restrict this classification to those torsion free
models where $\dim\{\mathfrak{K}\}>2$ to obtain models $\mathcal{M}_i^j(\star;0)$ where there
is an auxiliary parameter $\star$ in certain examples. If $j=6$, then $\dim\{\mathfrak{K}(\mathcal{M}_i^j(\star;0))\}=6$
and if $j=4$, then $\dim\{\mathfrak{K}(\mathcal{M}_i^j(\star;0))\}=4$. We then
add torsion to obtain models $\mathcal{M}_i^j(\star;T)$; we no longer have, of course,
that $\dim\{\mathfrak{K}(\mathcal{M}_i^j(\star;T))\}=j$ if $T\ne0$. Still, it seemed useful to keep the notation since
${}^0\!\mathcal{M}_i^j(\star;T)=\mathcal{M}_i^j(\star;0)$. We have that $\mathcal{M}_i^j(\star;T)$
and $\mathcal{M}_k^\ell(\star;\tilde T)$ are not linearly equivalent for $(i,j)\ne(k,\ell)$. Within a given class
defined by $(i,j)$ determining the precise set of representatives under linear equivalence is considerably more
delicate and we have not attempted such a finer classification.

We now establish the main result of the paper.
To obtain Assertion (1) in Theorem~\ref{T^1.3} we will compute the Lie algebras of affine Killing vector fields for all the models $\mathcal{M}_i^j(\star;T)$. We first write down a basis
for $\mathfrak{K}(\mathcal{M}_i^j(\star;0))$
and then examines the effect of the torsion tensor on the affine Killing equations
to derive the following result.

\begin{lemma}\label{L3.2}
	Let $\mathcal{M}$ be a Type~$\mathcal{A}$ model with torsion tensor $T=(T^1,T^2)$ so that
	$\dim\{\mathfrak{K}(\mathcal{M})\}>2$. Then $\mathcal{M}$ is linearly equivalent to one of the following surfaces with
	the values of $T$ listed; $\mathfrak{K}(\mathcal{M})=\mathrm{span}\{\partial_1,\partial_2\}$ for other values of $T$.%
	\begin{enumerate}\setlength\itemsep{0em}
		\item Let $\mathcal{M}_0^6(T):=\mathcal{M}(0, 0, T^1, T^2, -T^1, -T^2, 0, 0)$. Then
		\begin{enumerate}\setlength\itemsep{0em}
			\item $\mathfrak{K}(\mathcal{M}_0^6(0))=\mathrm{span}\{\partial_1,\partial_2,x^1\partial_1,x^1\partial_2,x^2\partial_1,x^2\partial_2\}$.
			\item $\mathfrak{K}(\mathcal{M}_0^6(T))=
			\mathrm{span}\{\partial_1,\partial_2, x^1(T^{1}\partial_1+T^{2}\partial_2), x^2(T^{1}\partial_1+T^{2}\partial_2)$ if $T\ne0$.
		\end{enumerate}
		\item Let $\mathcal{M}_1^6(T):=\mathcal{M}(1, 0, T^1, 1 + T^2, -T^1, 1 - T^2, 0, 0)$. Then  $T^1=0$ and
		\begin{enumerate}\setlength\itemsep{0em}
			\item $\mathfrak{K}(\mathcal{M}_1^6(0,0))=\mathrm{span}\{\partial_1,e^{-x^1}(\partial_1-x^2\partial_2),x^2\partial_2,
			{x^2(\partial_1-x^2\partial_2)},\partial_2,e^{-x^1}\partial_2\}$.
			\item 
			$\mathfrak{K}(\mathcal{M}_1^6(0,T^2))=
			\mathrm{span}\{\partial_1,\partial_2, x^2\partial_2, e^{-x^1}\partial_2\}$ if $T^2\ne0$.
		\end{enumerate}
		\item Let $\mathcal{M}_2^6(T):=\mathcal{M}({{{-1, 0}, {T^1, T^2}}, {{-T^1, -T^2}, {0, 1}}})$. Then $T^1T^2=0$ and
		\begin{enumerate}\setlength\itemsep{0em}
			\item $\mathfrak{K}(\mathcal{M}_2^6(0))=\mathrm{span}\{\partial_1,e^{x^1}\partial_1,e^{x^1+x^2}\partial_1,\partial_2,
			e^{-x^1-x^2}\partial_2,e^{-x^2}\partial_2\}$.
			\item $\mathfrak{K}(\mathcal{M}_2^6(T))=\mathrm{span}\{\partial_1,\partial_2,e^{-x^1-x^2}\partial_2,e^{-x^2}\partial_2\}$ if
			$T^1=0,T^2\ne0$.
			\item $\mathfrak{K}(\mathcal{M}_2^6(T))=\mathrm{span}\{\partial_1,e^{x^1}\partial_1,e^{x^1+x^2}\partial_1,\partial_2\}$ if
			$T^1\ne0,T^2=0$.
		\end{enumerate}
		\item Let $\mathcal{M}_3^6(T):=\mathcal{M}({{{0, 0}, {T^1, T^2}}, {{-T^1, -T^2}, {0, 1}}})$. 
		Then $T^1T^2=0$ and
		\begin{enumerate}\setlength\itemsep{0em}
			\item $\mathfrak{K}(\mathcal{M}_3^6(0))=\mathrm{span}\{\partial_1,x^1\partial_1,e^{x^2}\partial_1,\partial_2,e^{-x^2}\partial_2,
			x^1e^{-x^2}\partial_2\}$.
			\item $\mathfrak{K}(\mathcal{M}_3^6(T))=\mathrm{span}\{\partial_1,\partial_2,e^{-x^2}\partial_2,x^1e^{-x^2}\partial_2\}$
			if $T^1=0,T^2\ne0$.
			\item$\mathfrak{K}(\mathcal{M}_3^6(T))=\mathrm{span}\{\partial_1,\partial_2,x^1\partial_1,e^{x^2}\partial_1\}$
			if $T^1\ne0,T^2=0$.
		\end{enumerate}
		\item Let $\mathcal{M}_4^6(T):=\mathcal{M}({{{0, 0}, {T^1, T^2}}, {{-T^1, -T^2}, {1, 0}}} )$. Then $T^2=0$ and
		\begin{enumerate}\setlength\itemsep{0em}
			\item $\mathfrak{K}(\mathcal{M}_4^6(0))=\mathrm{span}\{\partial_1,\partial_2,(x^1+\frac12(x^2)^2)\partial_1,x^2\partial_1,\\
			\phantom{\mathfrak{K}(\mathcal{M}_4^6(0))=\mathrm{span}\{} ({ -}x^1x^2-\frac12(x^2)^3)\partial_1+(x^1+\frac12(x^2)^2)\partial_2,-(x^2)^2\partial_1+x^2\partial_2\}$.
			\item$\mathfrak{K}(\mathcal{M}_4^6(T^1,0))=\mathrm{span}\{\partial_1,\partial_2,(x^1+\frac12(x^2)^2)\partial_1,x^2\partial_1\}$ for $T^1\ne0$.
		\end{enumerate}
		\item Let $\mathcal{M}_5^6(T)=\mathcal{M}( {{{1, 0}, {T^1, 1 + T^2}}, {{-T^1, 1 - T^2}, {-1, 0}}})$. Then $T=0$ and\\
		$\mathfrak{K}(\mathcal{M}_5^6(0))=
		\mathrm{span}\{\partial_1,\partial_2,\cos(2x^2)\partial_1-\sin(2x^2)\partial_2,
		\sin(2x^2)\partial_1+\cos(2x^2)\partial_2$\\
		$\phantom{\mathfrak{K}(\mathcal{M}_5^6(0))=\mathrm{span}\{} e^{-x^1}(\cos(x^2)\partial_1-\sin(x^2)\partial_2),
		e^{-x^1}(\sin(x^2)\partial_1+\cos(x^2)\partial_2)\}$.
		\item Let $\mathcal{M}_1^4(T):=\mathcal{M}(-1, 0, 1 + T^1, T^2, 1 - T^1, -T^2, 0, 2)$. Then $T^2=0$ and\\
		$\mathfrak{K}(\mathcal{M}_1^4(T^1,0)) =
		\mathrm{span}\{\partial_1,\partial_2, e^{x^1}\partial_1, x^2e^{x^1}\partial_1\}$.
		\item Let $\mathcal{M}^4_2(c;T):=\mathcal{M}(-1, 0,c + T^1, T^2,c - T^1, -T^2,0, 1 + 2 c)$ for $c\ne0,-1$. Then $T^2=0$ and
		$\mathfrak{K}(\mathcal{M}^4_2(c;(T^1,0)))=
		\mathrm{span}\{\partial_1,\partial_2, e^{x^1}\partial_1, e^{x^1+x^2}\partial_1\}$.
		\item Let $\mathcal{M}^4_3(c;T):=\mathcal{M}(0, 0,c + T^1, T^2,c - T^1, -T^2,0, 1 + 2 c)$ for $c\ne0,-1$. Then $T^2=0$ and
		$\mathfrak{K}(\mathcal{M}^4_3(c;(T^1,0)))=
		\mathrm{span}\{\partial_1,\partial_2, e^{x^2}\partial_1, x^1\partial_1\}$.
		\item Let $\mathcal{M}^4_4(c;T):=\mathcal{M}(0, 0,1 + T^1, T^2,1 - T^1, -T^2,c, 2)$. Then $T^2=0$ and\\
		$\mathfrak{K}(\mathcal{M}^4_4(c;(T^1,0))) = \mathrm{span}\{\partial_1,\partial_2, \left(x^1+\tfrac{c}2(x^2)^2\right)\partial_1, x^2\partial_1\}$.
		\item Let $\mathcal{M}^4_5(c;T):=\mathcal{M}(1, 0,T^1, T^2,-T^1, -T^2,1 + c^2, 2 c)$. Then $T^2=0$ and\\
		$\mathfrak{K}(\mathcal{M}^4_5(c;(T^1,0))) = 
		\mathrm{span}\{\partial_1,\partial_2, e^{-x^1+c x^2}\cos{x^2}\partial_1,e^{-x^1+cx^2}\sin{x^2}\partial_1\}$.
	\end{enumerate}
\end{lemma}

One now performs a careful examination of the Lie algebras of Lemma~\ref{L3.2} to determine their isomorphism
type. This leads to the following classification result from which Theorem~\ref{T^1.3}~(1) follows:

\begin{lemma}
	Adopt the notation established in Lemma \ref{L3.2}. Let $\mathcal{M}$ be a Type~$\mathcal{A}$ model with torsion.
	Generically, $\mathfrak{K}(\mathcal{M})=\mathfrak{K}_{\mathcal{A}}$. 
	Let $\varepsilon\ne0$
	and let $(\varepsilon_1,\varepsilon_2)\ne(0,0)$.
	If $\dim\{\mathfrak{K}(\mathcal{M})\}>2$, then $\mathfrak{K}(\mathcal{M})$ has one of the following structures.
	\begin{enumerate}\setlength\itemsep{0em}
		\item $\mathfrak{K}(\mathcal{M}_i^6(0,0))\approx\mathbb{A}_6$.
		\item $\mathfrak{K}_{\mathcal{B}}\oplus\mathfrak{K}_{\mathcal{B}}\approx\mathfrak{K}(\mathcal{M}^6_1(0,\varepsilon)) \approx \mathfrak{K}(\mathcal{M}^6_2(0,\varepsilon)) \approx \mathfrak{K}(\mathcal{M}^6_2(\varepsilon,0)) \approx \mathfrak{K}(\mathcal{M}^6_3(\varepsilon,0))$\\
		$\phantom{\mathfrak{K}_{\mathcal{B}}\oplus\mathfrak{K}_{\mathcal{B}}}\approx\mathfrak{K}(\mathcal{M}^4_2(c;(T^1,0))) \approx \mathfrak{K}(\mathcal{M}^4_3(c;(T^1,0)))$.
		\item $A_{4,9}^0\approx\mathfrak{K}(\mathcal{M}_0^6(\varepsilon_1,\varepsilon_2))
		\approx\mathfrak{K}(\mathcal{M}^6_3(0,\varepsilon))\approx\mathfrak{K}(\mathcal{M}^6_4(\varepsilon,0))$\\
		$\phantom{A_{4,9}^0}\approx\mathfrak{K}(\mathcal{M}^4_1(c;(T^1,0)))\approx\mathfrak{K}(\mathcal{M}^4_4(c;(T^1,0)))
		$.
		\item $ A_{4,12}\approx\mathfrak{K}(\mathcal{M}^4_5(c;(T^{1},0)))$.
	\end{enumerate}
\end{lemma}

%%%%%%%%%%%%%%%%%%%%%%%%%%%%%%%%%%%%%%%%%%%%%%%%%%

\section{Type~$\mathcal{B}$ surfaces with torsion}\label{S4}
We proceed in a similar fashion in the Type~$\mathcal{B}$ setting. We parametrize the set of Type~$\mathcal{B}$
models by setting $\mathcal{N}(\vec\xi):=(\mathbb{R}^+\times\mathbb{R},\nabla^{\mathcal{B}}(\vec\xi))$ where
the Christoffel symbols take the form:
$$\begin{array}{llll}
\Gamma_{11}{}^1=(x^1)^{-1}\xi_1,&\Gamma_{11}{}^2=(x^1)^{-1}\xi_2,&\Gamma_{12}{}^1=(x^1)^{-1}\xi_3,&\Gamma_{12}{}^2=(x^1)^{-1}\xi_4,\\
\Gamma_{21}{}^1=(x^1)^{-1}\xi_5,&\Gamma_{21}{}^2=(x^1)^{-1}\xi_6,&\Gamma_{22}{}^1=(x^1)^{-1}\xi_7,&\Gamma_{22}{}^2=(x^1)^{-1}\xi_8.
\end{array}$$
The structure group for the set of Type~$\mathcal{B}$ models is not the full general linear group but rather the
$ax+b$ group acting by the shear $(x^1,x^2)\rightarrow(x^1,bx^1+ax^2)$; again we say two Type~$\mathcal{B}$ models
are {\it linearly equivalent} if they are in the same orbit of the induced linear action on $\mathbb{R}^8$.
The work of \cite{BVGRG18,GVR} mentioned
previously does not provide a full classification of all the Type~$\mathcal{B}$ models without torsion up to linear
equivalence. It does suffice, for our purposes, in that it does classify the torsion free Type~$\mathcal{B}$ models
with $\dim\{\mathfrak{K}\}>2$ by providing models $\mathcal{N}_i^j(\star;0)$ where $\star$ is an auxiliary parameter
present in some instances. Of particular interest are the geometries $\mathcal{N}_3^3$, which is the Lorentzian hyperbolic
plane, and $\mathcal{N}_4^3$, which is the hyperbolic plane.
The geometries $\mathcal{N}_i^4(\star;0)$ are also Type~$\mathcal{A}$
geometries. The torsion tensors are, of course, quite different. The geometries $\mathcal{N}_i^6(\star;0)$ are flat.
The proof of Lemma~\ref{L4.1} now follows by first writing down a basis for $\mathfrak{K}(\mathcal{N}_i^j(\star;0))$
and then examining the effect of the torsion tensor on the affine Killing equations.

\begin{lemma}\label{L4.1}
	Let $X:=x^1\partial_1+x^2\partial_2$.
	Let $\mathcal{N}$ be a Type~$\mathcal{B}$ model with torsion tensor $T=(T^1,T^2)$ so that
	$\dim\{\mathfrak{K}(\mathcal{N})\}>2$.  Then $\mathcal{N}$ is linearly equivalent to one of the following surfaces with
	the values of $T$ listed; $\mathfrak{K}(\mathcal{N})=\mathrm{span}\{X,\partial_2\}$ for other
	values of $T$.
	\begin{enumerate}\setlength\itemsep{0em}
		\item Let $\mathcal{N}_0^6(T):=\mathcal{N}({{{0, 0}, {T^1, T^2}}, {{-T^1, -T^2}, {0, 0}}})$. Then
		\begin{enumerate}\setlength\itemsep{0em}
			\item $\mathfrak{K}(\mathcal{N}_0^6(0))=\mathrm{span}
			\{\partial_1,\partial_2,x^1\partial_1,x^1\partial_2,x^2\partial_1,x^2\partial_2\}$.
			\item $\mathfrak{K}(\mathcal{N}_0^6(0,T^2))=
			\mathrm{span}\{X,\partial_2,x^1\partial_2,x^2\partial_2\}$ if $T^2\ne0$.
		\end{enumerate}
		\item Let $\mathcal{N}_1^6(\pm;T):=\mathcal{N}({{{1, 0}, {T^1, T^2}}, {{-T^1, -T^2}, {\pm1, 0}}})$.
		Then $T=0$ and\\
		\noindent$\mathcal{K}(\mathcal{N}_1^6(\pm,0))=\mathrm{span}\{X,\partial_2,\frac1{x^1}\partial_1,\frac{x^2}{x^1}\partial_1,
		\frac{(x^1)^2\pm(x^2)^2}{x^1}\partial_1,\frac{-x^2((x^2)^2\pm(x^1)^2)}{x^1}\partial_1+((x^1)^2\pm(x^2)^2)\partial_2\}$.
		\item Let $\mathcal{N}_2^6(c;T):=\mathcal{N}({{{-1 + c, 0}, {T^1, c + T^2}}, {{-T^1, c - T^2}, {0, 0}}})$ for $c\ne0$. Then
		$T^1T^2=0$ and
		\begin{enumerate}\setlength\itemsep{0em}
			\item $\mathfrak{K}(\mathcal{N}_2^6(c;0))=\mathrm{span}\{x^1\partial_1,x^2\partial_2,\partial_2,(x^1)^{-c}\partial_2,(x^1)^{-c}(x^1\partial_1
			-cx^2\partial_2), x^2(x^1\partial_1-cx^2\partial_2)\}$.
			\item $\mathfrak{K}(\mathcal{N}_2^6(c;(0,T^2)))= \mathrm{span}\{X,\partial_2, x^2\partial_2, (x^1)^{-c}\partial_2\}$ if $T^2\ne0$.
			\item $\mathfrak{K}(\mathcal{N}^6_2(-\tfrac12;{(T^1,0)})) = \mathrm{span}\{X,\partial_2, x^2(x^1\partial_1+\tfrac12 x^2\partial_2)\}$ if $T^1\ne0$.
		\end{enumerate}
		\item Let $\mathcal{N}_3^6(T):=\mathcal{N}({{{-2, 1}, {T^1, -1 + T^2}}, {{-T^1, -1 - T^2}, {0, 0}}})$. Then $T^1=0$ and
		\begin{enumerate}\setlength\itemsep{0em}
			\item $\mathfrak{K}(\mathcal{N}_3^6(0))=\mathrm{span}\{X,\partial_2,x^1\partial_2, (x^2+x^1\log x^1)\partial_2, \\
			\phantom{\mathfrak{K}(\mathcal{N}_3^6(0))=\mathrm{span}\{}{-}(x^1)^2\partial_1+ x^1(x^1-x^2)\partial_2,(x^2+x^1\log x^1)(-x^1\partial_1 + (x^1-x^2)\partial_2)\}$.
			\item $\mathfrak{K}(\mathcal{N}^6_3(0,T^2)) = \mathrm{span}\{X,\partial_2, x^1\partial_2, (x^2+x^1\log x^1)\partial_2\}$ if $T^2\ne0$.
		\end{enumerate}
		\item Let $\mathcal{N}_4^6(T):=\mathcal{N}({{{0, 1}, {T^1, T^2}}, {{-T^1, -T^2}, {0, 0}}})$. Then $T^1=0$ and
		\begin{enumerate}\setlength\itemsep{0em}
			\item $\mathfrak{K}(\mathcal{N}^6_4(0)) =\mathrm{span}\{X,\partial_2,x^1\partial_2, (x^2+x^1\log x^1)\partial_2, \\\phantom{\mathfrak{K}(\mathcal{N}^6_4(0)) =\mathrm{span}\{} \partial_1-(1+\log x^1)\partial_2,(x^2+x^1\log x^1)(\partial_1 {-} (1+\log x^1)\partial_2)\}$.
			\item $\mathfrak{K}(\mathcal{N}^6_4(0,T^2)) = \mathrm{span}\{X,\partial_2,x^1\partial_2, (x^2+x^1\log x^1)\partial_2\}$ if $T^2\neq0$.
		\end{enumerate}
		\item Let $\mathcal{N}_5^6(T):=\mathcal{N}({{{-1, 0}, {T^1, T^2}}, {{-T^1, -T^2}, {0, 0}}})$. Then $T^1=0$ and
		\begin{enumerate}\setlength\itemsep{0em}
			\item $\mathfrak{K}(\mathcal{N}^6_5(0)) =\mathrm{span}\{X,\partial_2,x^2\partial_2, x^1 x^2\partial_1,\log x^1\partial_2,x^1\log x^1\partial_{ 1}\}$.
			\item $\mathfrak{K}(\mathcal{N}^6_5(0,T^2)) = \mathrm{span}\{X,\partial_2,x^2\partial_2,\log x^1\partial_2\}$ if $T^2\neq0$.
		\end{enumerate}
		\item Let $\mathcal{N}_6^6(c;T):=\mathcal{N}({{{c, 0}, {T^1, T^2}}, {{-T^1, -T^2}, {0, 0}}})$ for $c\ne0,-1$. Then $T^1=0$ and
		\begin{enumerate}\setlength\itemsep{0em}
			\item $\mathfrak{K}(\mathcal{N}^6_6(c;0)) =\mathrm{span}\{X,\partial_2,x^2\partial_2, (x^1)^{-c}\partial_1, (x^1)^{-c}x^2\partial_1, (x^1)^{c+1}\partial_2\}$.
			\item $\mathfrak{K}(\mathcal{N}^6_6(c;(0,T^2))) = \mathrm{span}\{X,\partial_2,x^2\partial_2,(x^1)^{c+1}\partial_2\}$ if $T^2\neq0$.
		\end{enumerate}
		\item Let $\mathcal{N}_1^4(\kappa;T):=
		\mathcal{N}({{{2 \kappa, 1}, {T^1, T^2 + \kappa}}, {{-T^1, -T^2 + \kappa}, {0,  0}}})$
		for $\kappa\ne0,-1$. Then $T^1=0$ and
		$\mathfrak{K}(\mathcal{N}_1^4(\kappa;(0,T^2)))=
		\mathrm{span}\{X,\partial_2, x^1\partial_2, x^1(\partial_1-\log x^1 \partial_2)\}$.
		\item Let $\mathcal{N}_2^4(\kappa,\theta;T):=\mathcal{N}({{{-1 + \theta + 2 \kappa, 0}, {T^1, 
					T^2 + \kappa}}, {{-T^1, -T^2 + \kappa}, {0, 0}}})$ for $\theta\ne0$ and $\kappa(\kappa+\theta)\ne0$. Then $T^1=0$ and
		$\mathfrak{K}(\mathcal{N}_2^4(\kappa,\theta;(0,T^2)))
		= \mathrm{span}\{X,\partial_2, x^2\partial_2, (x^1)^{\theta}\partial_2\}$.
		\item Let $\mathcal{N}_3^4(c;T):=\mathcal{N}({{{2c-1, 0}, {T^1, 
					T^2 + c}}, {{-T^1, -T^2 +c}, {0, 0}}})$ for $c\ne0$. Then $T^1=0$ and
		$\mathfrak{K}(\mathcal{N}_3^4(c;(0,T^2)))=
		\mathrm{span}\{X,\partial_2, x^2\partial_2, \log x^1\partial_2\}$.
		\item Let $\mathcal{N}_1^3(\pm;T):=\mathcal{N}(-\frac32, 0, T^1, -\frac12 + T^2, -T^1, -\frac12 - T^2, \pm\frac12, 0)$.
		Then $T^2=0$ and $\mathfrak{K}(\mathcal{N}_1^3(\pm;(T^1,0)))=
		\mathrm{span}\{X,\partial_2, x^2\left(2x^1\partial_1+x^2\partial_2\right)\}$.
		\item Let $\mathcal{N}_2^3(c;T):=
		\mathcal{N}({{{-\frac32, 0}, {1 + T^1, -\frac12 + T^2}}, {{1 - T^1, -\frac12 - T^2}, {c, 2}}})$. Then $T^2=0$ and
		$\mathfrak{K}(\mathcal{N}_2^3(c;(T^1,0)))
		= \mathrm{span}\{X,\partial_2,  x^2\left(2x^1\partial_1+x^2\partial_2\right)\}$.
		\item
		Let $\mathcal{N}_3^3(T):=\mathcal{N}({{{-1, 0}, {T^1, -1 + T^2}}, {{-T^1, -1 - T^2}, {-1, 0}}})$. Then $T=0$ and\\
		$\mathfrak{K}(\mathcal{N}_3^3(T))=\mathrm{span}\{X,\partial_2,
		2x^1x^2\partial_1+((x^2)^2+(x^1)^2)\partial_2\}$.
		\item Let $\mathcal{N}_4^3(T):=\mathcal{N}({{{-1, 0}, {T^1, -1 + T^2}}, {{-T^1, -1 - T^2}, {1, 0}}})$. Then $T=0$
		and\\ $\mathfrak{K}(\mathcal{N}_4^3(0))=\mathrm{span}\{X,\partial_2,
		2x^1x^2\partial_1+((x^2)^2-(x^1)^2)\partial_2\}$.
	\end{enumerate}
\end{lemma}

One now performs a careful examination of the Lie algebras of Lemma~\ref{L4.1} to determine their isomorphism
type. This leads to the following classification result from which Theorem~\ref{T^1.3}~(2) follows:

\begin{lemma}
	Adopt the notation established in Lemma~\ref{L4.1}. Let $\mathcal{N}$ be a Type~$\mathcal{B}$ model
	with torsion. Generically, $\mathfrak{K}(\mathcal{N})$  is 2-dimensional and is isomorphic to the 2-dimensional non-Abelian
	Lie algebra $\mathfrak{K}_{\mathcal{B}}$. Let $\varepsilon\ne0$. 
	If $\dim\{\mathfrak{K}(\mathcal{N})\}>2$, then $\mathfrak{K}(\mathcal{N})$ has one of the following structures.
	\begin{enumerate}\setlength\itemsep{0em}
		\item $\mathfrak{K}(\mathcal{N}_i^6(0,0))\approx\mathbb{A}_6$.
		\item $\mathfrak{K}_{\mathcal{B}}\oplus\mathfrak{K}_{\mathcal{B}}\approx\mathfrak{K}(\mathcal{N}^6_0((0,\varepsilon)))\approx\mathfrak{K}(\mathcal{N}^6_2( c;(0,\varepsilon)))\approx \mathfrak{K}(\mathcal{N}_3^6(0,\varepsilon))\approx \mathfrak{K}(\mathcal{N}_4^6(0,\varepsilon))$ \\ 
		$\phantom{\mathfrak{K}_{\mathcal{B}}\oplus\mathfrak{K}_{\mathcal{B}}} \approx \mathfrak{K}(\mathcal{N}_6^6(0,\varepsilon))\approx \mathfrak{K}(\mathcal{N}^4_1(c;(0,T^2)))\approx \mathfrak{K}(\mathcal{N}^4_2(\kappa,\theta; (0,T^2)))$.
		\item $ A_{4,9}^0\approx\mathfrak{K}(\mathcal{N}^6_5(0,\varepsilon))
		\approx \mathfrak{K}(\mathcal{N}^4_3(c;(0,T^2)))$.
		\item $\mathfrak{so}(2,1)\approx\mathfrak{K}(\mathcal{N}^6_2(-\tfrac12; (\varepsilon,0))) \approx \mathfrak{K}(\mathcal{N}^3_1(\pm; (T^1,0))) \approx \mathfrak{K}(\mathcal{N}^3_2(c;(T^1,0)))$ \\ 
		$\phantom{\mathfrak{so}(2,1)} \approx \mathfrak{K}(\mathcal{N}^3_3(0,0)) \approx \mathfrak{K}(\mathcal{N}^3_4(0,0))$.
	\end{enumerate}
\end{lemma}

Note that for the particularly interesting cases $\mathcal{N}^3_3$ and $\mathcal{N}^3_4$ (the Lorentzian and Riemannian hyperbolic planes) 
any torsion perturbation reduces their Lie algebra of affine Killing vectors from $\mathfrak{so}(2,1)$ to $\mathfrak{K}_{\mathcal{B}}$. These two surfaces, 
together with $\mathcal{N}^6_1(\pm)$, are the only cases of homogeneous Type $\mathcal B$ surfaces which under any perturbation with a torsion 
tensor reduces the Lie algebra of affine Killing vectors to $\mathfrak{K}_{\mathcal{B}}$.

%%%%%%%%%%%%%%%%%%%%%%%%%%%%%%%%%%%%%%%%%%%%%%%%%%

\section{Conclusions}

Possible extensions of General Relativity are based on the independence between the metric and the affine properties of spacetime. 
In this context torsion plays a fundamental role. In the present article we examine the effects of torsion on the affine Killing vectors of a surface. 
Since we consider homogeneous surfaces we have a large number of symmetries that preserve the affine connection. In fact, even flat surfaces with non-zero torsion tensor have a very rich structure.

In this paper we have obtained a complete description of the Lie algebra $\mathfrak{K}(\mathcal M)$ of affine Killing vectors fields on any homogeneous surface $\mathcal M$ with non-vanishing torsion. In the Type $\mathcal A$ setting $\mathfrak{K}(\mathcal M)$ is restricted to be one of the following: $\mathfrak{K}_{\mathcal B}\oplus\mathfrak{K}_{\mathcal B}$, $A^0_{4,9}$, $A^4_{12}$, or  $\mathfrak{K}_{\mathcal A}$. In the Type $\mathcal B$ setting, $\mathfrak{K}(\mathcal M)$ can only be one of the following: $\mathfrak{K}_{\mathcal B}\oplus\mathfrak{K}_{\mathcal B}$, $A^0_{4,9}$, $\mathfrak{so}(2,1)$, or $\mathfrak{K}_{\mathcal B}$. This completes the analysis of \cite{BVGRG18}.

We believe that a systematic description of affine structures with torsion is useful in the search of interesting 
non-metrizable geometries. A detailed classification of homogeneous surfaces in terms of the torsion tensors 
they admit is currently in progress. 
There is no immediate extension of this work to the higher-dimensional setting since there is no analogous
classification of the possible 
affine models, even if torsion is absent. However, we recall that Lemma \ref{L1.4} holds in any dimension; 
it is plausible that the analysis of the addition of torsion to a given torsion-free connection at the level of the affine Killing equations 
could give some insight on possible approaches to the problem.

\begin{acknowledgments}
Research of DD was partially supported by Universidad Nacional de La Plata under grant 874/18 and project 11/X791. Research of PBG was partially supported by Project MTM2016-75897-P (AEI/FEDER, Spain). Research of PP was partially supported by a Fulbright-CONICET scholarship and by Universidad Nacional de La Plata under project 11/X615. DD and PP thank the warm hospitality at the Mathematics Department of the University of Oregon, where this work was carried out.
\end{acknowledgments}

%%%%%%%%%%%%%%%%%%%%%%%%%%%%%%%%%%%%%%%%%%%%%%%%%%

%%%%%%%%%%%%%%%%%%%%%%%%%%%%%%%%%%%%%%%%%%%%%%%%%%

\end{document}